\documentclass[leqno,12pt]{article}

\RequirePackage[OT1]{fontenc}
\RequirePackage{amsthm,amsmath}
\RequirePackage[colorlinks,citecolor=blue,urlcolor=blue]{hyperref}

\usepackage{latexsym}
\usepackage{amsmath}
\usepackage{amssymb}
\usepackage{bbm}
\usepackage{color}
\usepackage{graphicx}
\usepackage{epsfig}
\usepackage[latin1]{inputenc}

\title{Large deviations and Wschebor's theorems}

\author{\small{Jos\'e R. Le\'on}\footnote{IMERL, Universidad de la Rep\'ublica, Montevideo, Uruguay and Escuela de Matem\'atica, Universidad
Central de Venezuela, rlramos@fing.edu.uy}
\and {\small Alain Rouault}\footnote{Laboratoire de Math{\' e}matiques de Versailles, UVSQ, CNRS, Université Paris-Saclay, 78035-Versailles Cedex France, alain.rouault@uvsq.fr}}

\begin{document}
\maketitle
\begin{abstract}
We revisit Wschebor's theorems on small increments for processes with scaling and stationary properties and deduce large deviation principles.
\end{abstract}
\bigskip

\smallskip

{\bf MSC 2010:} {0J65, 60G51, 60F15, 60F10}

\newcommand\Z{{\mathbb Z}}
\newcommand\N{{\mathbb N}}
\newcommand\R{{\mathbb R}}
\newcommand\C{{\mathbb C}}
\newcommand\T{{\mathbb T}}
\renewcommand\P{{\mathbb P}}
\newcommand{\demi}{\frac{1}{2}}
\newcommand\ind{{\mathbf 1}}
\newcommand\E{{\mathbb E}}
\newcommand\eps{\varepsilon}
\newcommand{\Cov}{\mbox{\rm Cov}}
\newcommand{\Var}{\mbox{\rm Var}}
\newcommand{\Corr}{\mbox{\rm Corr}}
\newcommand{\ve}{\varepsilon}

\def \el{\sur{=}{(d)}}
\def \be{\begin{eqnarray*}}
\def \ee{\end{eqnarray*}}
\def \ben{\begin{eqnarray}}
\def \een{\end{eqnarray}}
\def\wh{\widehat}
\def\wt{\widetilde}
\def\pn{\par\noindent}
\def\sn{^{(N)}}
\def\tr{\hbox{tr}\!\ }
\def\AArm{\fam0 }
\def\AAk#1#2{\setbox\AAbo=\hbox{#2}\AAdi=\wd\AAbo\kern#1\AAdi{}}%
\def\AAr#1#2#3{\setbox\AAbo=\hbox{#2}\AAdi=\ht\AAbo\raise#1\AAdi\hbox{#3}
}%
\def \eref#1{(\ref{#1})}
\def \sur#1#2{\mathrel{\mathop{\kern 0pt#1}\limits^{#2}}}
\def\proof{\noindent{\bf Proof:}\hskip10pt}
\def\QED{\hfill\vrule height 1.5ex width 1.4ex depth -.1ex \vskip20pt}
\newcommand{\ii}{{\mathrm{i}}}
\newcommand{\cCn}{\mathcal{C}_n}
\newcommand{\br}{\mathcal{B}}
\newcommand{\usn}{a_{n}}
\newcommand{\Int}{\mathrm{int}\,}
\newcommand{\Clo}{\mathrm{clo}\,}
\newcommand{\prs}[2]{\langle#1,#2\rangle}
\def\cro#1{\llbracket #1  \rrbracket}

\newtheorem{lem}{Lemma}[section]
\newtheorem{defi}[lem]{Definition}
\newtheorem{theo}[lem]{Theorem}
\newtheorem{cor}[lem]{Corollary}
\newtheorem{prop}[lem]{Proposition}
\newtheorem{nota}[lem]{Notation}
\newtheorem{rem}[lem]{Remark}
\newtheorem{rems}[lem]{Remarks}
\newtheorem{conj}[lem]{Conjecture}
\def\proof{\noindent{\bf Proof:}\hskip10pt}
\setcounter{page}{1}
\setcounter{section}{0}

\normalsize\rm

\section{Introduction :   Wschebor's theorem and beyond}
In 1992, Mario Wschebor \cite{wschebor1992accroissements} proved the following  remarkable property of the linear Brownian motion 
 $(W(t) , t\ge0 ; W(0)=0)$.  If 
 $\lambda$ is the Lebesgue measure on $[0,1]$, then,  
 almost surely, for every $x \in \Bbb R$ and 
 every $t \in [0,1]$:
\begin{eqnarray}\label{W95}\lim_{\ve \to 0} \lambda \{s\le t:\, \frac{W(s+\ve) -W(s)}{\sqrt\ve}\le x\} =  t\Phi(x)\,,  
\end{eqnarray} 
where $\Phi(x) = \mathcal N((-\infty, x])$ and $\mathcal N$ is the standard normal distribution. 
It  is a sort of law of large numbers (LLN) for the random measure defined as 
\[\mu_{\varepsilon}(A)=\lambda\{ s\in[0,1]:\,\frac{W(s+\varepsilon)-W(s)}{\sqrt\varepsilon}\in A\}\,,\]
 that a.s. weakly converges towards $\mathcal N$.

This result was generalized 
shortly after by Wshebor \cite{wschebor1995almost} for L\'evy processes and by Aza\"{\i}s and Wschebor \cite{azais1996almost} for  
random processes with stationary increments and other processes.  Moreover, also in \cite{wschebor1992accroissements} and \cite{wschebor1995almost},  the result was shown for mollified processes as follows. 

For $\psi\in BV$ (the set of bounded variation functions on $\mathbb R$),  
 let  
\[\psi^\varepsilon (t) = \frac{1}{\varepsilon} \psi\left(\frac{t}{\varepsilon}\right)\]
denote the rescaled version of $\psi$  
and for $X$ a measurable function, set $X_\psi^\varepsilon = X \star \psi^\varepsilon$, i.e.
\begin{eqnarray}
\label{star1}X_\psi^\varepsilon (t) := \int \psi^\varepsilon (t-s) X(s) ds\,,\end{eqnarray}
and
\begin{eqnarray}\label{star2}\dot X_\psi^\ve :=\int X(s) d\psi^\ve (t-s)\,.\end{eqnarray}

The result  reads, when $X = W$ and $\psi \in BV \cap L^2$,
\begin{eqnarray}\label{asconverg}\lambda\{s\le t:\, \mathcal W_\psi^\ve (s) 
\le x\}\to t\Phi(x/ ||\psi||_2)\ \ \ (a.s.)\,,\end{eqnarray} 
with 
\begin{eqnarray}
\label{star3}\mathcal W_\psi^\ve (s) := \sqrt\varepsilon\!\ \dot W_\psi^{\varepsilon}(s)\,.\end{eqnarray}

Notice that when $\psi = \psi_1 :=1_{[-1, 0]}$, 
\[\dot X_\psi^\ve(t) = \ve^{-1}\left(X(t+\ve) - X(t)\right)\ , \ \mathcal W_\psi^\ve := \ve^{-1/2} (W(\cdot + \ve) - W(\cdot))\,,\]
and we recover (\ref{W95}). 

In subsequent articles the a.s. result was extended to obtain a stable central limit theorem  (CLT). For instance, let us consider a real even function $g$ such that  $\mathbb E[g^2(N)]<\infty,$ for $N \sim \mathcal N$ 
 (a typical example is $g(x)=|x|^p-\mathbb E[|N|^p]$). Defining for a bounded and continuous function $f$ the family 
$$\mathcal Z_{\varepsilon}=\frac1{\sqrt\ve}\int_0^1f(W(s))\left(g\left(\mathcal W_\psi^\ve (s)
\right)-\mathbb Eg(N)\right)ds,$$
we observe first that the convergence in (\ref{asconverg}) implies that
$$\lim_\ve \int_0^1f(W(s))g\left(\mathcal W_\psi^\ve (s)
\right)ds =  \mathbb Eg(N/ ||\psi||_2) \int_0^1f(W(s))ds \ \ \ (a.s.)\,.$$
Moreover, 
for $\mathcal S$ the convergence stable of measures we have 
\begin{eqnarray}\label{BL}\lim_\ve \mathcal Z_{\varepsilon}\stackrel{\mathcal S}=  \sigma(g)\int_0^1 f(W(s))dB(s),\end{eqnarray}
where $B(s)$ is another Brownian motion independent of $W$ and $\sigma$ is an explicit positive constant. Let us point out that if we take $f=1$ and by integrating on the interval $[0,t]$ the above result 
turns into  a functional CLT: 
\begin{eqnarray}\label{CW}\lim_\ve \left( \frac1{\sqrt\ve}\int_0^t\left(g\left(\mathcal W_\psi^\ve (s)
\right)-\mathbb Eg(N)\right)ds, t \in [0,1]\right)=  \left(\sigma(g) B(t), t \in [0,1]\right)\,,\end{eqnarray}
in distribution.

The result in (\ref{BL}) was obtained in \cite{berzin1997brownian}. Since then,  such type of matters were generalized to: diffusions, fractional Brownian motion (fBm), stationary increments Gaussian processes, L\'evy processes, etc. A very complete review 
with a large number of references can be found in \cite{wschebor2006smoothing}. More recently,  
 in 2008, Marcus and Rosen in \cite{Marcus2008SPA} have studied the convergence of the $L^p$  norm  (this is $g(x)=|x|^p$ 
 in (\ref{CW})) of the increments of stationary  Gaussian processes. In the cited article the authors closed the problem in a somewhat definitive form. 
In another article (\cite{Marcus2008AP}) they  said that their proofs were initially based on  Wschebor's method, but afterwards they changed and looking for a more general and broadly used procedure.

When we are faced with a LLN-type result (a.k.a. convergence of a family of random objects to a deterministic one), it is nowadays natural to ask for a possible large deviation principle (LDP). 

Let us give some notations. If $\Sigma= \mathbb R, \mathbb R^+ \times \mathbb R$ or $[0,1]\times \mathbb R$, we denote by $\mathcal M^+(\Sigma)$ and $\mathcal M^r (\Sigma)$ the set of  Borel measures on $\Sigma$ positive and having total mass $r$, respectively.

If $Z$ is a measurable function from $\mathbb R^+$ to $\mathbb R$, let $M_Z \in \mathcal M^+(\mathbb R^+ \times \mathbb R)$ be defined by
\begin{eqnarray}
\label{defmu1}
M_Z( I \times A)  =\lambda\{ s\in I:\, Z(s)
\in A\}\,,
\end{eqnarray}
for every Borel subset $I \times A$ of $\mathbb R^+ \times \mathbb R$.
The first marginal of $M_Z$ is $\lambda$. 
The second marginal $\mu_Z$  is defined either by its action on a Borel set $A$
\begin{eqnarray}
\label{defmarg}
\mu_Z(A) = M_Z ([0,1] \times A)=   \lambda\{s \in [0,1] : Z \in A\}
\end{eqnarray}
or, by its action on a test function $f \in \mathcal C_b(\mathbb R)$ (set of  bounded continuous functions on $\mathbb R$) 
$$\int_{\mathbb R} f(x)  d\mu_Z (x) = \int_0 ^1 f(Z(t)) \ dt  
\,,$$
so that $\mu_{Z}$ is the occupation measure
\[\mu_{Z} =   \int_0 ^1 \delta_{Z} dt\,.\]

In this framework, 
we can consider (\ref{W95})  and (\ref{asconverg}) as laws of large numbers (LLN):
\[M_{\mathcal W_1^\ve} \Rightarrow \lambda \times \mathcal N\ \ (a.s.)\]
where $\Rightarrow$ stands for the weak convergence. 
Since the Brownian motion $W$ is self-similar (Property P1) and has stationary increments (P2), it is possible to reduce the problem about $\mu_{\mathcal W_1^\varepsilon}$ ($\varepsilon \rightarrow 0$)
 to a problem of an occupation measure in large time ($T :=\varepsilon^{-1} \rightarrow \infty$) for a process $Y$ independent of $\varepsilon$.  This new process is stationary and ergodic. Moreover the independence of increments of $W$ (P3) and its self-similarity induces a 1-dependence for $Y$, which allows to apply a criterion of Chiyonobu and Kusuoka \cite{chiyonobu1988large} to get an LDP.

Actually, as the crucial properties (P1, P2, P3) are shared by $\alpha$-stable L\'evy processes, we state the LDP in this last framework. This is the content of Section 3 with an extension to random measures built with mollifyers.  Previously a basic lemma on equalities in law is stated in Section 2.

The fBM  with Hurst index $H \not= 1/2$ shares also properties (P1, P2) but not (P3) with the above processes. Nevertheless, since it is Gaussian, with an explicit spectral density, we prove the LDP for $(\mu_\varepsilon)$ 
under specific conditions on the mollifier, thanks to a criterion of \cite{Bryc-Dembo}. This is the content of Section 4.  In Section \ref{SpaT} we state an LDP for the space-time measure  defined in (\ref{defmu1}) when $Z$ is one of the above processes and in Section \ref{LevP}, we state a result for some ``process level" empirical measure.
At last, 
in Section \ref{DisV} we  study  discrete versions of Wschebor's theorem.

Among the issues not addressed here, we may quote:  increments for Gaussian random fields in $\mathbb R^d$ and  
multi-parameter indexed processes.

Let us notice that except in a specific case in Section \ref{Looking}, we cannot give an explicit expression for the rate function. Moreover
if one would be able to prove that the rate function is strictly convex and its minimum is reached at $\lambda \times \mathcal N$, this would give an alternate proof of  Wschebor's results.

\section{General framework}

Recall that a real-valued process $\{X(t), t \in \mathbb R\}$ 
\begin{itemize}
\item
has stationary increments if
\[\{X(t+h) - X(h) , t\in \mathbb R\} \el \{X(t) - X(0) , t\in \mathbb R\}\,,\]
\item
 is self-similar with index $H > 0$ if
 \[\forall a > 0 \ \{X(at) , t \in \mathbb R\} \el \{a^H X(t) , t \in \mathbb R\}\,.\]
\end{itemize}
If $X$ is a self-similar process with index $H$ we set,
 if $\psi \in BV$
\begin{eqnarray}
\mathcal X_\psi^\ve= \ve^{1-H} \dot X_\psi^\ve\,,
\end{eqnarray}
where $\dot X_\psi^\ve$ is defined as in (\ref{star2}) by
\begin{eqnarray}
\dot X_\psi^\ve= \int X(s) d\psi^\ve(t-s) = \frac{1}{\ve}\int X(t- \ve u) d\psi(u)\,.
\end{eqnarray}
 In particular
\begin{eqnarray}\label{defpsi1}
\mathcal X_\psi^1(t) = \int X(s) d\psi(t-s)\,.
\end{eqnarray}
The following lemma is the key for our study. 
\begin{lem}
\label{cruz}
Assume that $X$ is self-similar with index $H$.
For fixed $\ve$ and $\psi \in BV$, we have
\begin{eqnarray}
\label{inlaw}\left(\mathcal X_\psi^\ve(t) , t \in \mathbb R\right) &&\el \left(\mathcal X_\psi^1(t\ve^{-1}) , t \in \mathbb R\right)\\
\mu_{\mathcal X_\psi^\ve} &&\el \ve \int_0^{\ve^{-1}}  \delta_{\mathcal X_\psi^1 (t)} dt\,.
\end{eqnarray}
Moreover, if $X$ has stationary increments, then 
$\mathcal X_\psi^1$ is stationary.
\end{lem}
\proof 
It  is straightforward. First,
\begin{eqnarray}
\dot X_\psi^\ve (t)&= \ve^{-1} \displaystyle\int X(t- \ve u) d\psi(u) \el \ve^{H -1} \displaystyle\int X\left(\frac{t}{\ve} - u\right) d\psi(u)\,,
\end{eqnarray}
where the last equality comes from self-similarity and holds as a process in $t \in \mathbb R$. This yields (\ref{inlaw}),
 and then
\begin{eqnarray}\mu_{\mathcal X_\psi^\ve} &= \displaystyle\int_0^1 \delta_{\mathcal X_\psi^\ve(t)} dt= 
\ve  \displaystyle\int_0^{1/\ve}  \delta_{\mathcal X_\psi^\ve(\ve \tau)} d\tau
\el \ve \displaystyle\int_0^{1/\ve} \delta_{\mathcal X_\psi^1(\tau)} d\tau\,. \ \ \qed
\end{eqnarray}
We give now a definition which will set the framework for 
 the processes studied in the sequel. Recall that the $\tau$-topology on $\mathcal M^1(\mathbb R)$  is the topology   induced by the space of bounded measurable functions on $\mathbb R$. It is stronger than the weak topology which is induced by $\mathcal C_b (\mathbb R)$.

\begin{defi}
Let $\mathcal F$ be a subset of the set $BV$ of bounded variation function from $\mathbb R$ in $\mathbb R$. We say that a process $X$ with stationary increments and self-similar with index $H$ has the  $(LDP_w,  \mathcal F, H)$ (resp. $(LDP_\tau,  \mathcal F, H)$) property if the process $\mathcal X_\psi^1$ is well defined and if  for every $\psi \in \mathcal F$,  the family $(\mu_{\mathcal X_\psi^\ve})$ satisfies the LDP in $\mathcal M^1(\mathbb R)$ equipped with the weak topology (resp. the $\tau$-topology), in the scale $\ve^{-1}$,
 with good rate function 
\begin{eqnarray}\label{g15'}\Lambda_\psi^* (\mu) = \sup_{f \in \mathcal C_b(\mathbb R)} \int f d\mu - \Lambda_\psi (f)\,,
\end{eqnarray}
(the Legendre dual of $\Lambda_\psi$) where for $f\in \mathcal C_b(\mathbb R)$,
\begin{eqnarray}
\label{g15''}
\Lambda_\psi(f) = \lim_{T\rightarrow \infty} T^{-1}\log \mathbb E \exp \int_0^T f(\mathcal X_\psi^1(t)) dt\,,\end{eqnarray}
in particular, the above limit exists.
\end{defi}
Roughly speaking, this means that the probability of seeing $\mu_{\mathcal X_\psi^\ve}$ close to $\mu$ for a small $\ve$ is asymptotically $e^{-\Lambda_\psi^*(\mu)/\ve}$.
\section{The $\alpha$-stable L\'evy process}
Let $\alpha \in (0,2]$ fixed. The $\alpha$-stable L\'evy process $(S(t) , t\geq 0; S(0) = 0)$ has independent and stationary increments and is $1/\alpha$-self-similar. 
If $\psi\in BV$ is compactly supported, we set
\[
\mathcal S_\psi^\ve(t) := 
\varepsilon^{1-1/\alpha} \int S(s) d\psi_\ve(t-s)\,\]
and as in (\ref{defmu1}-\ref{defmarg}), we may build the measures $M_{\mathcal S_\psi^\ve}$
 and $\mu_{\mathcal S_\psi^\ve}$.
In \cite{azais1996almost}, Theorem 3.1, it is proved that a.s.
\[M_{\mathcal S_\psi^\ve} \Rightarrow \lambda \times \Sigma_\alpha\ (a.s.)\]
where $\Sigma_\alpha$ is the law of $||\psi||_\alpha S(1)$.
\begin{prop}
If $\mathcal F$ is the set of bounded variation functions with compact support, then 
the $\alpha$-stable L\' evy process has the $(LDP_\tau,  \mathcal F, 1/\alpha)$  property.
\end{prop}

\proof 
 We apply Lemma \ref{cruz} with $X = S$ and $H = 1/\alpha$.  

Assume that the support of $\psi$ is included in $[a,b]$. 
Since $S$ has independent and stationary increments, the process $\mathcal S_\psi^1$
is stationary and $(b-a)$-dependent. This last property means that $\sigma (\mathcal S_\psi^1(u) , u \in A)$ and $\sigma (\mathcal S_\psi^1(u), u \in B)$ are independent as soon as the distance between $A$ and $B$ is greater than $1$. 
The process $(\mathcal S_\psi^1)$ satisfies the condition (T) in 
  \cite{chiyonobu1988large} 
 (see also condition (S) in    
 \cite{Bryc-Dembo_h} p. 558). Then the family 
$(\mu_{\mathcal S_\psi^\ve})$ 
satisfies the LDP and the other conclusions hold. \qed
\begin{rem}
When $\alpha=2$ we recover the Brownian case. In particular, when $\psi =\psi_1$ 
\begin{eqnarray}\label{defSl}\mathcal S_\psi^1 (u)
= W(u+1) -W(u)\ , \ u \in\mathbb R\,.\end{eqnarray} 
This process is called often Slepian process; 
it  is Gaussian, stationary and
1-dependent.
\end{rem}

\section{The fractional Brownian motion}

\subsection{General statement}
We treat now the case of self-similar Gaussian processes with stationary increments, i.e. fractional Brownian motion (fBm in short). The fBm  with Hurst parameter $H \in [0,1)$ is the Gaussian process $(B_H(t), t \in \mathbb R)$ with covariance
\[\mathbb E B_H(t) B_H(s) = \frac{1}{2}\left(|s|^{2H}+|t|^{2H} - |t-s|^{2H}\right)\,,\]
It has a chaotic (or harmonizable) representation of $B_H$ (see \cite{sam} Prop. 7.2.8)
\begin{eqnarray}
\label{harmo}B_H (t) = \frac{1}{C_H} \int_\mathbb R \left(e^{i\lambda t} - 1\right) |\lambda|^{-H-\frac{1}{2}} d{\bf W}(\lambda)\end{eqnarray}
where $\bf W$ is a complex Brownian motion and 
\[C_H^2 = 
\frac{2\pi}{\Gamma(2H +1) \sin(\pi H)}\,.\]

This process has stationary increments and is self-similar of index $H$.
When $H=1/2$ we recover the Brownian motion, and it is the only case where the increments are independent.

When $\psi \in BV$ with compact support, the LLN can be formulated as :
\begin{eqnarray}
M_{\mathcal X_\psi^\ve}\Rightarrow \lambda \times  \mathcal N( \cdot \  \sigma_\psi)\ \ \ (a.s.)\,,  
\end{eqnarray}
where 
\[\sigma^2_\psi = - \frac{1}{2}\int\!\!\! \int  |u-v|^{2H} d\psi(u) d\psi(v)\,,\]
(see \cite{azais1996almost}).

Our result on large deviations is the following. In Fourier analysis  
 we adopt the following notation:
 when $f,g \in L^1(\mathbb R)$
\[\hat f(\theta) = \int e^{it\theta} f(t) dt \ , \ \check g (\gamma) = \frac{1}{2\pi} \int e^{-i\gamma x}g(x) dx\,.\] 

\begin{prop}
\label{4.1}
Denote
\[\mathcal G :=  \{ \psi \in BV \cap L^1 
\}\]
\[\mathcal G_H = \{ \psi\in L^1 : \exists \lim_{\lambda\rightarrow 0} |\hat\psi(\lambda)||\lambda|^{\frac{1}{2} -H}\}\ , \ (0< H < 1)\,.\] 
The process $B_H$ has the $(LDP_w, \mathcal F, H)$ property if one of the following conditions are satisfied:
\begin{eqnarray}
\label{H2}
H \leq 1/2 &&\hbox{and} \  \ \mathcal F = \mathcal G
\,,\\
\label{H3}
1/2< H < 1&&\hbox{and} \  \
\mathcal F = 
\mathcal G \cap \mathcal G_H
\,.
\end{eqnarray}
\end{prop}
Particular cases are examined in Section \ref{Particular}.
\begin{rem}
\label{4.2}
If we define
\[\mathcal G_0 = \{\psi \in L^1 :  \int \psi(t) dt = 0 \ \hbox{and} \ \int |t\psi(t)| dt < \infty
\}\]
then it holds that
\begin{eqnarray}
\label{GyG} \mathcal G_0 \subset \mathcal G_1 \subset\mathcal G_H\,.\end{eqnarray}
\end{rem}

\noindent\underbar{Proof of Proposition  \ref{4.1}}:

We apply Lemma \ref{cruz} with $X = B_H$. But now, for lack of independence, we will work with the spectral density and apply Theorem  2.1 in \cite{Bryc-Dembo}, which ensures the LDP as soon as the spectral density is in $\mathcal C_0(\mathbb R)$.

Using (\ref{harmo}), the process $\mathcal X_\psi^1$ may be written as 
\begin{eqnarray}
\mathcal X_\psi^1(t) = \int B_H(t-s)d\psi(s) = C_H^{-1} \int\!\!\!  \int  \left(e^{\ii \lambda(t- s)} - 1\right) |\lambda|^{-H-\frac{1}{2}} d{\bf W}(\lambda) d\psi(s)\,. 
\end{eqnarray}
Now, when $\psi \in \mathcal G$ it holds that  $ \lim_{|t| \rightarrow \infty} |\psi(t)| = 0$ and by integration by parts, 
\begin{eqnarray}
\label{IPP}
\int e^{i\lambda s} d\psi(s) &=& -\ii\lambda\hat\psi(\lambda)\\
\nonumber\int  \left(e^{\ii \lambda (t-s}) - 1\right)d\psi(s) &=& \ii \lambda e^{\ii \lambda t} \hat\psi(-\lambda)\end{eqnarray}
so that
\begin{eqnarray}
\mathcal X_\psi^1(t) = \ii  C_H^{-1}\int e^{\ii t \lambda} \hat\psi(-\lambda)\frac{\lambda}{|\lambda|^{H+\frac12}} d{\bf W}(\lambda)\,.
\end{eqnarray}
 The spectral density of the stationary process $\mathcal X_\psi^1$ is then
 \begin{eqnarray}
\label{SDfBm}\ell_H(\lambda)=C_H^{-2} \left| \hat\psi(\lambda)\right|^2
|\lambda|^{1- 2H}
\,.\end{eqnarray}
Let us check conditions on $H$ and $\psi$ so that  $\ell_H \in \mathcal C_0(\mathbb R)$.

\begin{itemize}
\item
From (\ref{IPP}) it holds that 
 $|\hat\psi(\lambda)| = 0 (|\lambda|^{-1})$ and then, for all $0 < H < 1$
\begin{eqnarray}\label{0infty}\lim_{|\lambda|\rightarrow \infty} |\ell_H(\lambda)| = 0\,.\end{eqnarray}
\item
If $\psi \in L^1$,  $\hat \psi\in \mathcal C_0$  (Riemann-Lebesgue) and $\ell_H$ is even and 
continuous on $(0,\infty)$.
\item

For $H \leq 1/2$, the continuity of $\ell_H$ at $0$ is obvious.

For $H > 1/2$, this continuity is ensured by the assumption $\psi \in \mathcal G_H$. \qed
\end{itemize}
\medskip

\noindent\underbar{Proof of Remark  \ref{4.2}}:

To prove (\ref{GyG}), note that $\mathcal G_1 \subset \mathcal G_H$ is obvious and that, 
 if $\psi \in \mathcal G_0$,
\begin{eqnarray*}\label{first}
\hat\psi(0)  &=& \int \psi(x) dx = 0\\
|\hat\psi(\lambda)| &=& \left|\int (e^{\ii\lambda x}-1) \psi(x) dx\right| \leq |\lambda|\int |x\psi(x)| dx\,,
\end{eqnarray*}
so that $\psi \in \mathcal G_1$. \qed

\subsection{Contraction}
Since the mapping $\mu \mapsto \int |x|^p d\mu(x)$ is not continuous, 
we cannot obtain an LDP  for the moments of $\mu_{\mathcal X_\psi^\ve}$ by invoking the contraction principle (Th. 4.2.1 in  \cite{DZ10})). Nevertheless, in the case of the fBm, the Gaussian stationary character of the process allows to conclude
. It is a direct application of Corollary 2.1 in \cite{Bryc-Dembo}.
\begin{prop}
If either $H \leq 1 /2$ and $\psi \in \mathcal G$ or $H > 1/2$ and $\psi \in  \mathcal G \cap \mathcal G_H$, then 
the family $\left(\int_0^1 |\mathcal X_\psi^\ve(t)|^2 dt \right)$, where $X=B_H$,  satisfies the LDP, in the  scale $\ve^{-1}$ with good rate function 
\[I_\psi (x) = \sup_{-\infty < y < 1/(4\pi M)} \{xy - L(y)\}\,,\]
where
\[L(y) = -\frac{1}{4\pi}\int \log (1-4\pi y \ell_H(s)) ds\]
and
\[M = \sup_\lambda \ell_H(\lambda)\,.\]
More generally, for $0 \leq p \leq 2$, the family $\left(\int_0^1 |\mathcal X_\psi^\ve(t)|^p dt \right)$ satisfies the LDP at scale $\ve$ with a convex rate function.
\end{prop}

\subsection{Particular cases}
\label{Particular}

\subsubsection{Two basic mollifiers}

1) As seen before, the function $\psi_1= 1_{[-1,0]}$ is the most popular. It  allows to study the first order increments $X(t+\ve) - X(t)$. It belongs to $\mathcal G$ but since
\[|\hat\psi_1(\lambda)| =\frac{|\sin (\lambda/2)|}{|\lambda/2|}\]
it does not belong to $\mathcal G_H$ for $H>1/2$.

 For $H = 1/2$, we recover the Brownian motion and replace the notation $\mathcal X$ by $\mathcal W$. The process $\mathcal W_{\psi_1}^1$ is the Slepian process (\ref{defSl}) 
 with covariance
\[r(t) = \left(1 - |t|\right)^+\,,\]
and 
spectral density: 
\[
  \check r (\lambda) = \frac{1}{2\pi}\left(\frac{\sin \frac{\lambda}{2}}{\frac{\lambda}{2}}\right)^2 
\,.\]
As it is said above since $\check r$ 
  is $\mathcal C_0$, 
the  occupation
  measure  satisfies a LDP in the weak topology in the scale $\ve^{-1}$.  in the scale $\varepsilon^{-1}$.
This argument could have been used to prove the LDP, instead of the argument in Section 2 (but for the weak topology and not the $\tau$-topology). Notice that although 
$\check r$  is differentiable, we could not apply Theorem 5.18 in Chiyonobu and Kusuoka \cite{chiyonobu1988large}, since the condition (5.19) therein is violated in $x \in 2\pi\mathbb Z$.

2) Another interesting function is
\[\psi_2 = \frac{1}{2}\left(1_{[-1, 0]}-1_{[0,1]}\right)\,\]
which yields
\begin{eqnarray}
\label{2nd}
\dot X_{\psi_2}^\ve(t) =  \frac{ X(t+\ve)-2X(t) + X(t-\ve)}{2\ve}\,.\end{eqnarray}
Since 
\[|\hat\psi_2(\lambda)| = \frac{\sin^2 (\lambda/2)}{|\lambda/2|}\,,\]
we see that $\psi_2 \in \mathcal G \cap \mathcal G_H$ for every $H \in (0,1)$ and then $(\mu_{\mathcal X_{\psi_2}^\ve})$ satisfies 
the LDP.

In (\ref{2nd}) we are faced with  second order increments of the process $X$.
 These increments are linked with the behavior of the second derivative of $X^{\eps}$ when it exists. Let us consider $\psi$ smooth enough  so that $X^\ve_\psi$, defined in (\ref{star1}), has a second derivative. 
 For instance, let $\psi\in  \mathcal G$ 
 and such that $\psi' \in \mathcal G$. 
Then  the function 
$X^\ve_\psi$  is twice differentiable and
\begin{eqnarray*}
\ddot X_\psi^\ve(t) = \ve^{-2} \int X(t-\ve s) d\psi'(s) = \ve^{-1} \dot X_{\psi'}^\ve (t)\,.  
\end{eqnarray*}
Now, $\psi' \in \mathcal G_H$ since
\[|\widehat \psi' (\lambda)| |\lambda|^{\frac{1}{2}-H} = |\hat \psi (\lambda)| |\lambda|^{\frac{3}{2}-H}\rightarrow 0\] as $\lambda \rightarrow 0$.

Since $\mathcal X_{\psi'}^\ve = \ve^{2-H}\ddot X^\ve_\psi$,  we conclude that for every $H \in (0,1)$, the family  $(\mu_{\ve^{2-H} \ddot X_\psi^\ve})$ satisfies the LDP in the  scale $\ve^{-1}$ and good rate function $\Lambda^\star_{\psi'}$. The choice 
\[\psi(t) = \frac{1}{2}\left(1 - |t|\right)^+\]
allows to recover $\psi' = \psi_2$ and the second order increments.

\subsubsection{Looking for an explicit rate function}
\label{Looking}

It is not easy to find examples of explicit rate functions for the occupation measures of the above stationary processes $\mathcal X_\psi^1$, since in general  the limiting cumulant generating function $\Lambda$ is not explicit. A particularly nice situation in the Gaussian case will occur if the process is also Markovian, i.e. if $\mathcal X_\psi^1$ is the Ornstein-Uhlenbeck (OU) process. Indeed, for the OU, the rate function for the LDP of the occupation measure is given by the Donsker-Varadhan theory (\cite{stroock2012introduction} ex. 8.28) :
$$\Lambda^{*}(\mu)=\frac{1}{2}\int_{\R}|g'(x)|^2 d\mathcal N(x)$$
if $d\mu = g^2 d\mathcal N$. The goal is then to find a mollifier $\psi$ such that $\mathcal X_\psi^1$ is distributed as OU. 

To begin with, let us assume that the underlying process is Brownian, which implies that $\mathcal W_\psi^1$
  is  again Gaussian and stationary, with
  spectral density (cf. (\ref{SDfBm})):
\begin{eqnarray*}
\label{F1}
\check r(\lambda) 
 = \frac{1}{2\pi} |\hat\psi(\lambda)|^2\,.
\end{eqnarray*}
\medskip

 For OU, the covariance and spectral density are, respectively
\[r(t) = e^{-|t|}\ , \ \check r(\lambda) = \frac{1}{\pi(1 + \lambda^2)}\,.\]
To solve the equation
\[\mathcal X^1_\psi \el OU\] 
turns out to solve
\begin{eqnarray}
\label{square root}
|\hat\psi(\lambda)|^2 =  \frac{2}{1+\lambda^2}\,.
\end{eqnarray}

We present two answers. 

1) Let us choose
\begin{eqnarray*}
\hat\psi(\lambda)= \frac{\sqrt 2}{1- i\lambda} \ , \ \psi (x) = 
\sqrt 2 e^{-x} 1_{[0, \infty)}(x)\,,
\end{eqnarray*}
and then, the formula (\ref{star3}) becomes
\[\mathcal W_\psi^1 (t) = \sqrt 2 \int_{-\infty}^t e^{-(t-s)} dW_s\]
which is the classical representation of the stationary OU as a stochastic integral (\cite{sam} p.138).

2) Let us choose $\psi$ such that
\begin{eqnarray*}
\hat\psi(\lambda) = \frac{\sqrt 2}{\sqrt{1+ \lambda^2}}
\end{eqnarray*}
This is equivalent to say
\[\psi(x) = \frac{1}{2\pi}\int  e^{-ix\lambda} \frac{\sqrt 2}{\sqrt{1+ \lambda^2}} d\lambda\]
i.e.
\begin{eqnarray*}\psi(x) = \frac{\sqrt 2}{\pi}\int _0^\infty \frac{\cos(x\lambda)}{\sqrt{1+\lambda^2}}  d\lambda
= \frac{\sqrt 2}{\pi} K_0(x),\end{eqnarray*}
where $K_0$  is 
 the MacDonald (or modified  Bessel) function (see \cite{Davies} p.369 or \cite{Fouriertransforms} formula 17 p.9). This function  can be expressed also as
\[K_0(x) = \sqrt \pi e^{-x} \Psi (1/2, 1 ;2x)\,,\]
where $\Psi$ is the confluent hypergeometric function (see \cite{Bateman} p. 265).
\medskip

Let us now extend the study to the fBm. Looking for a kernel $\psi$  leading to the OU process, (\ref{SDfBm}) leads to the equation
\begin{eqnarray*}
\label{LFOU}
\left| \hat\psi(\lambda)\right|^2 = C_H^{2} \frac{|\lambda|^{2H-1}}{\pi(1+ \lambda^2)}\,,
\end{eqnarray*}
hence, for instance if $\psi$ is even,
\begin{eqnarray*}
\label{raizLFOU}
 \hat\psi(\lambda) = C_H \frac{|\lambda|^{H-\frac12}}{\sqrt{\pi(1+\lambda^2)}}\,.
\end{eqnarray*}
For $H< 1/2$, we did not find a closed expression for the kernel 
\begin{eqnarray}
\label{LDPc}
\psi (x) = \frac{C_H}{\pi}\int_0^\infty \cos(\lambda x)\frac{|\lambda|^{H-\frac12}}{\sqrt{\pi(1+\lambda^2)}} d\lambda
\end{eqnarray}
  in the literature. 

For $H > 1/2$, this function is not continuous in $0$, so it cannot be the Fourier transform of an integrable kernel.
We have proved
\begin{prop}
When $H \leq 1/2$ and $\psi$ is given by (\ref{LDPc}),  the family $(\mu_{\mathcal H_\psi^\ve})$ satisfies the LDP, in the scale $\ve^{-1}$ with good rate function
$$\Lambda^{*}(\mu)=\frac{1}{2}\int_{\R}|g'(x)|^2 d\mathcal N(x)\,$$
if $d\mu = g^2 d\mathcal N$.
\end{prop}
\begin{rem} 
In this case, $\Lambda^*$ has a unique minimum at $\mu = \mathcal N$ which allows to recover Wschebor's result on a.s. convergence.
\end{rem}

\section{A space-time LDP} 
\label{SpaT}
We will state a complete LDP for some of our models, i.e. an LDP for $(M_{\mathcal X_\psi^\ve})$, whenever  $(\mu_{\mathcal X_\psi^\ve})$ satisfies the LDP. Following the notations of Dembo and Zajic in \cite{DZa} we denote by $\mathcal A\mathcal C_0$ the set of maps $\nu : [0,1]\rightarrow \mathcal M^+ (\mathbb R)$ 
such that
\begin{itemize}
\item $\nu$ is absolutely continuous with respect to the variation norm,
\item $\nu(0) = 0$ and  $\nu(t) - \nu(s) \in \mathcal M^{t-s}(\mathbb R)$ for all $t> s \geq 0$,
\item for almost everty $t\in [0,1]$, $\nu(t)$ possesses a weak derivative.
\end{itemize}
(This last point means that $\nu(t+\eta)- \nu(t)/\eta$ has a limit as $\eta \rightarrow 0$ - denoted by $\dot\nu(t)$- in $\mathcal M^+(\mathbb R)$ equipped with the topology of weak convergence).

 Let $F$ 
\begin{eqnarray}\nonumber\mathcal M^+([0,1]\times \Bbb R) &&\rightarrow D\left([0,1] ; \mathcal M^+(\Bbb R)\right)\\
M &&\mapsto \left(t \mapsto F(M)(t) = M([0,t] , \cdot)\right)
\end{eqnarray}
or in other words $F(M)(t)$ is the positive measure on $\Bbb R$ defined by its action on $\varphi \in \mathcal C_b$:
\[\langle F(M)(t) , \varphi\rangle = \langle M , 1_{[0,t]} \times \varphi\rangle\,.\]
Here $D([0,1] ; \cdot)$ is the set of c\`ad-l\`ag functions, equipped with the supremum norm topology. At last, let $\mathcal E$ be the image of $F$. 

\begin{theo}
When the process $X$ is the $\alpha$-stable L\' evy process and $\psi\in BV$ has a compact support, the family $\left(M_{\mathcal X_\psi^\ve}\right)$ satisfies the LDP in $\mathcal M^1([0,1]\times\mathbb R)$ equipped with the weak topology, in the scale $\ve^{-1}$ with the good rate function
\begin{eqnarray}
\label{ratem}
\Lambda^*(M)  & = \begin{cases}  \displaystyle\int_0^1 \Lambda_\psi^* (\dot\gamma(t)) dt& \text{ if }  \gamma := F(M) \in \mathcal A\mathcal C_0, \\
\infty & \text{ otherwise. } 
\end{cases}
\end{eqnarray}
\end{theo}

\proof
As in the above sections, it is actually a problem of large deviations in large time. 
For the sake of simplicity, set
\[Y= \mathcal X_\psi^1\]
and $T= \ve^{-1}$.
Using Lemma \ref{cruz}, the problem reduces to the study of
the family $(M_{Y(\cdot T)})$. First, we  study the corresponding distribution functions.

Actually, we have
\begin{eqnarray}\label{Micmac}F(M_{Y(\cdot T)})(t) 
 = \int_0^t \delta_{Y(sT)} ds= T^{-1}\int_0^{tT} \delta_{Y(s)} ds =:H_T (t) \,.\end{eqnarray}
In a first step we will prove that the family $(H_T)$ 
 satisfies the LDP, then it a second step we will transfer this property to $M_{Y(\cdot T)}$.
\smallskip

\noindent\underline{First step} : We follow the method of Dembo-Zajic \cite{DZa}. We begin with a 
reduction to their ``discrete time" method by introducing
\[\xi_k = \int_{k-1}^k \delta_{Y_s} ds , (k \geq 1)\ \ \hbox{and} \ \ S_T (t) = \sum_1^{\lfloor tT\rfloor} \xi_k\,.\]
It holds that
\begin{eqnarray} T^{-1}\int_0^{tT} \delta_{Y(s)}ds - T^{-1} S_T (t) 
 = T^{-1}\int_{\lfloor tT\rfloor}^{tT} \delta_{Y(s)} ds\end{eqnarray}
and this difference has a total variation norm less than $T^{-1}$, so that the families $(T^{-1}S_T)$ and $(H_T)$
 are exponentially equivalent (Def. 4.2.10 in \cite{DZ10}). 

The sequence $\xi_k$ is 1-dependent, hence satisfies 
 condition (S) in \cite{DZa} p.22 
which implies, by Th.  4 in the same paper 
 that $(T^{-1}S_T)$ satisfies the LDP in $D([0,1]; \mathcal M^+(\mathbb R))$ provided with the uniform norm topology,  with the convex good rate function 
\begin{eqnarray}
\label{Igrande}
I(\nu) = \int_0^1 \Lambda_\psi^\star (\dot\nu(t)) dt\end{eqnarray}
when $\nu \in \mathcal A\mathcal C_0$ and $\infty$ otherwise.
 
We conclude, owing to Th. 4.2.13 in \cite{DZ10}, that $(H_T)$ satisfies the same LDP. 
\medskip

\noindent\underline{Second step} :
 We have now to carry this LDP to $(M_{Y(\cdot T)})$ (see (\ref{Micmac})). 
For every $T >0$, $H_T \in \mathcal E
 \subset D([0,1) ; \mathcal M^+(\mathbb R)$. 
We saw that the effective domain of $I$ is included in $\mathcal E$. 
So, by Lemma 4.1.5 in Dembo-Zeitouni \cite{DZ10}, $(H_T)$ satisfies the same LDP in $E$
 equipped with the (uniform) induced topology. Now, $F$ is bijective from $\mathcal M^1([0,1] \times \mathbb R)$ to $\mathcal E$. Let us prove that $F^{-1}$ is continuous from $\mathcal E$ (equipped with the uniform topology) to $\mathcal M^1([0,1] \times \mathbb R)$ equipped with the weak topology. 

For $f : [0,1] \rightarrow \mathbb R$, let 
\begin{eqnarray}
\Vert f \Vert_{BL} &&= \sup_x |f(x)| + \sup_{x \not= y} \frac{|f(x) - f(y)|}{|x-y|}\\
\label{dBL}
d_{BL}( \mu, \nu) &&= \sup_{f : \vert f\Vert_{BL} \leq 1} \left|\int f d\mu - \int fd\nu\right|
\end{eqnarray}
The space $\mathcal M^+(\mathbb R)$ is a Polish space when equipped with  the topology induced by $d_{BL}$, compatible with the weak topology.

It is known that $M_n \rightarrow M\in \mathcal M^1([0,1] \times \mathbb R)$ weakly as soon as
\begin{eqnarray}M_n(1_{[0,t]} \otimes f) \rightarrow M(1_{[0,t]} \otimes f)
\end{eqnarray}
for every $t \in [0,1]$ and every $f$ such that $\Vert f \Vert_{BL} < \infty$. 
But, for such $t,f$ 
we have
\begin{eqnarray}
\sup_t |M_n(1_{[0,t]} \otimes f) - M(1_{[0,t]} \otimes f)|\leq  d_{BL}(F(M_n), F(M))
\end{eqnarray}
which implies that $F^{-1}$ is continuous from $\mathcal E$ to $\mathcal M^1([0,1] \times \mathbb R)$.

By the contraction principle (Th. 4.2.1 in \cite{DZ10}) we deduce that $M_{Y(\cdot T)}$ satisfies the LDP in $\mathcal M^1([0,1] \times \mathbb R)$ with good rate function $J(M) = I(F(M))$, wherer $I$ is given by (\ref{Igrande}). \qed
\medskip

\section{``Level process" study}
\label{LevP}
In the study of strong convergence problems such as the a.s. CLT (see \cite{Heck} and \cite{March}), an interesting problem is the LDP of empirical measures at the level of processes.  If we restrict us to the Brownian case to simplify, the corresponding problem  could be the behavior of 
\[\int_0^1 \delta_{\left\{\frac{W(s+ \ve) - W(s)}{\sqrt \ve}, \ s \geq t\right\}} dt\,.\]
Here we do not see clearly the interest of such a study for the Wschebor's theorem. 
Nevertheless, it seems natural to consider  the family $(\xi^\ve_t,  t\geq 0)$ of shifted processes 
\begin{eqnarray}
\label{defLve}
\xi^\ve_t : s \mapsto \frac{W(t+\ve s) - W(t)}{\sqrt \ve} \in \mathcal C([0,1]) \,,
\end{eqnarray}and 
the following empirical measure
\begin{eqnarray}
\mathcal L_\ve :=  \int_0^1 \delta_{\xi^\ve_t}\!\  dt\,.
\end{eqnarray}
By the scaling invariance, for every $\ve >0$, 
\begin{eqnarray}
\label{scale}
(\xi^\ve_{\ve t} , t \geq 0) \el (\xi^1_t , t \geq 0)\,,
\end{eqnarray}
and then
\begin{eqnarray}\label{eqlaw}\mathcal L_\ve = \int_0^1 \delta_{\xi^\ve_t}\!\  dt\el \tilde{\mathcal L}_\ve := \ve \int_0^{\ve^{-1}} \delta_{\xi_t^1
} dt\,.\end{eqnarray}
Since we have
\begin{eqnarray}
\label{mS}
\xi_t^1 = (W(t+s) - W(t) , s\in [0,1]) ,
\end{eqnarray}
 the process $(\xi^1_{t} , t \geq 0)$ will be called the 
 the {\it meta-Slepian} process in the sequel. For every $t$, the distribution of $\xi_t^1$ is  the Wiener measure $\mathbb W$ on $\mathcal C([0,1])$.


The meta-Slepian process is clearly stationary and 1-dependent. 
Since it  is ergodic,  the Birkhoff theorem tells us that, almost surely when $\ve\to 0$, $\tilde{\mathcal L}_\ve$ converges weakly to $\mathbb W$.
From the equality in distribution (\ref{eqlaw}) we deduce that $(\mathcal L_\ve)$ converges in distribution to the same limit. But this limit is deterministic, hence the convergence of  $(\mathcal L_\ve)$ holds in probability.
We just proved:
\begin{theo}
When $\ve \to 0$, the family of random probability measures $(\mathcal L_\ve)$ on $\mathcal C([0,1])$ converges  in probability weakly to the Wiener measure $\Bbb W$ on $\mathcal C([0,1])$.
\end{theo}
The problem of almost sure convergence raises some difficulties. We have obtained on the one hand  a partial almost sure fidi convergence (which is no more that a multimiensional extension of Wschebor theorem) and on the other hand an almost sure convergence when we plug $\mathcal C([0,1])$ into the Hilbert space $L^2([0,1])$, equipped with its norm.  

To this last purpose, 
if $\mu$ is a measure on $\mathcal C([0,1])$, we will denote by $\mu^L$ its extension to $L^2([0,1)]$, i.e. that for every $B$ Borel set of  $L^2([0,1])$,
\[\mu^L(B) = \mu(B \cap \mathcal C([0,1])\,.\]
\begin{theo}
\label{cvproc}
\begin{enumerate}
\item 
For every integer $d$ and every $t_1, \dots, t_d \in [0,1]$, almost surely when $\ve \to 0$, the family $(\mathcal L_\ve\pi_{t_1, \dots, t_d}^{-1})$ of random probability measures on $\Bbb R^d$ converges  weakly to $\Bbb W\pi_{t_1, \dots, t_d}^{-1}$ on $\mathcal C([0,1])$, where $\pi_{t_1, \dots, t_d}$ be the projection : $f \in \mathcal C([0,1]) \mapsto f(t_1), \dots, f(t_d)$.
\item 
When $\ve \to 0$, the family of random probability measures $(\mathcal L_\ve^L)$ on $L^2([0,1])$ converges  weakly almost surely to the Wiener measure $\Bbb W^L$ on $L^2([0,1])$.
\end{enumerate}
\end{theo}

We failed to prove a (full) almost sure fidi convergence, i.e. in 1. to state that ``almost surely, for every $t_1, \dots, t_d$ ...". Moreover we do not know if an almost sure convergenge at the level of processes is true. 

For the proof, we need the following lemma, which is straightforward owing to the properties of stationarity and 1-dependence.
\begin{lem}
\label{az}
If $F$ is a bounded differentiable function with bounded derivative from $\mathcal C([0,1])$ 
(resp. $L^2([0,1])$) to $\mathbb R$. Then
\begin{eqnarray}
a.s. \ \lim_{\ve \rightarrow 0}\int_0^1 F(\xi_t^\ve)\!\ dt = \int_{\mathcal C([0,1])} F(\xi) \mathbb W(d\xi)\,.
\end{eqnarray}
\end{lem}

\noindent\underline{Proof of Lemma \ref{az}}:

It is along the lines  of \cite{azais1996almost}. We first claim a quadratic convergence as follows. By Fubini and stationarity
\[\mathbb E\left(\int_0^1 F(\xi_t^\ve) dt\right) = \int_0^1 \mathbb E  F(\xi_t^\ve) dt = \int_{\mathcal C([0,1])} F(\xi) \mathbb W(d\xi)\,,\]
and by Fubini and $1$-dependence,
\begin{eqnarray}\Var \left(\int_0^1 F(\xi_t^\ve) dt\right) = \int\!\int_{|t-s| < 2\ve} \Cov\left(F(\xi_t^\ve), F(\xi_s^\ve)\right) dt ds \leq 4\ve ||F||_\infty
\,.\end{eqnarray}
The Borel-Cantelli lemma implies a.s. convergence of $\int_0^1 F(\xi_t^\ve) dt$ along 
any sequence $(\ve_n)$ such that $\sum_n \ve_n < \infty$.

To go on, take $\ve_{n+1} < \ve < \ve_n$ and notice that
\begin{eqnarray}
\label{83}
\left|\int_0^1 F(\xi^\ve_t) - F(\xi_t^{\ve_n})\right| dt \leq ||F'||_\infty\sup_{t,u} \left|\xi_t^\ve(u) - \xi_t^{\ve_n}(u)\right|\,. 
\end{eqnarray}
Now we use the properties of Brownian paths. 
On the interval $[0, 2]$ the Brownian motion satisfies a.s. a Holder condition with exponent $\beta< 1/2$, so that we can define the a.s. finite random variable 
\begin{eqnarray}
M := 2\sup_{u,v \in [0,2]}\frac{|W(u) - W(v)|}{|v-u|^\beta }
\,.
\end{eqnarray}
So, 
\begin{eqnarray}
\nonumber
\sup_{s\in [0,1]}|\xi_t^\ve(s) - \xi_t^{\ve_n} (s)|&&\displaystyle\leq \frac{M}{2} \frac{(\ve_n - \ve)^\beta}{\ve^{1/2}} +  \frac{M}{2}(\ve_n)^\beta \left(\ve^{-1/2} - (\ve_n)^{-1/2}\right)\\\nonumber
&&\displaystyle= \frac{M}{2}\frac{(\ve_n)^\beta}{\ve^{1/2}}\left[\left(1 - \frac{\ve}{\ve_n}\right)^\beta + \left(1 - \sqrt{\frac{\ve}{\ve_n}}\right)\right]
\\\label{magicW}&&\leq M \frac{\ve_n^\beta - \ve^\beta}{\ve^{1/2}} \leq  M \frac{\ve_n^\beta - \ve_{n+1}^\beta}{\ve_{n+1}^{1/2}}\,.
\end{eqnarray}
The choice of  $\ve_n = n^{-a}$ with $a>1$ and $\beta \in \left(\frac{a}{2(a+1)}, \frac{1}{2}\right)$ ensures that 
the right hand side of (\ref{magicW}), hence of (\ref{83}) tends to $0$ a.s., which ends the proof. \qed
\bigskip

\noindent\underline{Proof of Theorem \ref{cvproc}}

1. The (random) characteristic functional of the (random) probability measure $\mathcal L_\ve$ on $[0,1]$ equipped with the Borel $\sigma$-field and the Lebesgue meausre is a function from the dual space of 
$\mathcal C([0,1])$, i.e. $\mathcal M([0,1])$ to $\mathbb C$ defined by
\[\hat G_\ve :  \rho  \mapsto \int_0^1\exp \{\ii \int \xi^\ve _t(s) \rho(ds)\} dt\,.\]
Actually,  $\hat G_\ve(\rho) = \int F(\xi_t^\ve) dt$ with $F(\xi) = \exp \{\ii \int \xi(u) \rho(du)\}$. This function fulfils the conditions of Lemma \ref{az}. 

We have then, for every $\rho$, 
\begin{eqnarray}a.s. \ \lim \hat G_\ve(\rho) = \hat G(\rho)&&:= \int_{\mathcal C([0,1])} \exp \{\ii \int_0^1 \xi(s) \rho(ds)\} \mathbb W(d\xi)\\ &&= \exp\left(-\frac{1}{2}\int_0^1 (\rho([u, 1])^2 du\right).\end{eqnarray}

Let fix $d$ and $t_1, \dots, t_d$. For every  ${\bf a} := (a_1, \dots ,a_d)\in \Bbb R^d$, let us consider the 
measure $\rho_{{\bf a}} =\sum_1^d a_k \delta_{t_k}$ and the
following event
\[A({\bf a}) := \left\{\lim_{\ve\rightarrow 0}
\hat G_\ve(\rho_{{\bf a}}) 
 = \hat G \left(\rho_{{\bf a}}\right)\right\}.\]
The above analysis tells us that $\Bbb P(A({\bf a}))= 1$ for every $\bf a$. 
By a classical argument using Fubini's theorem we deduce that 
almost surely, for almost every $ {\bf a} \in \Bbb R^d$
\begin{eqnarray}
  \lim_{\ve\rightarrow 0}\int_0^1 \exp\left(\ii \sum_1^d a_k \xi_t^\ve(t_k)\right) dt = \hat G \left(\sum_1^d a_k \delta_{t_k}\right).
\end{eqnarray}  
By a slight adaptation of the L\' evy's continuity theorem (which is detailed in Appendix), we conclude that $(\mathcal L_\ve\pi_{t_1, \dots, t_d}^{-1})$ converges weakly to the good limit.
\medskip

2. We will use a method coming from \cite{gine1980central} p. 46\footnote{It is used there to prove that in Hilbert spaces, convergence in the Zolotarev metric implies weak convergence.} . It consists in checking Billingsley's criterion on intersection of balls (\cite{billingsleyprobability} p.18) and approximating indicators by smooth functions. Let us give details for only one ball to shorten the proof.

For $\delta \in (0,1)$, define
\begin{eqnarray}\phi_{\delta}(t)=\mathbf 1_{(-\infty,1]}(t)+\mathbf 1_{[1,(1+\delta)^2]}(t) \frac1C\int_0^{\frac{((1+\delta)^2 -t)}{(2\delta+\delta^2)}}e^{ -\frac1{s(1-s)}}ds\,\,,\end{eqnarray} where
$$C=\int_0^1e^{-\frac1{s(1-s)}}ds\,.$$
The function $\phi_{\delta}$ has a bounded support and it is continuous and $||\phi_{\delta}||_{\infty}=1$. Now we consider
$\psi_\delta :  L^2([0,1]\to\R$ defined by
$$\psi_\delta(\xi)=\phi_{\delta}(||\xi||^2).$$This function is $\mathcal C^\infty$ and 
 has all its derivatives bounded. For every $\xi_c \in L^2([0,1]), r > 0, \delta \in (0,r)$ we have the nesting
\begin{eqnarray}
\mathbf 1_{B(\xi_c; r-\delta)} (\xi)\leq \psi_{\frac{\delta}{r-\delta}}\left(\frac{\xi-\xi_c}{r-\delta}\right) \leq \mathbf 1_{B(\xi_c;r)}(\xi)
\leq \psi_{\frac{\delta}{r}}\left(\frac{\xi-\xi_c}{r}\right) \leq \mathbf 1_{B(\xi_c; r+\delta)}(\xi)\,.\nonumber\\
\label{nesting}
\end{eqnarray}
Take a sequence $\delta_n \to 0$.

Let us remind that the measure $\mathcal L_\ve^L$ is random. We did not write explicitly the item $W$  for simplicity, although it is present in (\ref{defLve}).

For every test function $F$ as in Lemma \ref{az}, we have a null set $N_F$ such that for $W \notin N_F$ 
\begin{eqnarray}
\label{Nunico}\int_{L^2([0,1])} F(\xi) \mathcal L_\ve^L(d\xi)  \to \int_{\mathcal C([0,1])} F(\xi) \mathbb W (d\xi)\,.\end{eqnarray}
Let $(g_k)_{k\geq 1}$ be a countable dense set in $L^2([0,1])$, and for $q \in \mathbb Q$, 
\[F^-_{n,k,q}(\xi) = \psi_{\delta_n/(q-\delta_n)}\left(\frac{\xi-g_k}{q-\delta_n}\right) , F^+_{n,k,q}(\xi) =\psi_{\delta_n/q}\left(\frac{\xi-g_k}{q}\right)\]
and
\[N =\bigcup_{n,k,q} \left(N_{F^-_{n,k,q}}\cup N_{F^+_{n,k,q}} \right)\,. \]
Take $W \notin N$. Assume that the ball $B(\xi_c ; r)$ is given. Take $\gamma >0$, then by density one can find  $k\geq 1$ and $q\in \mathbb Q^+$ such that
\begin{eqnarray}
\label{ballq}||\xi_c -g_k|| \leq \gamma \ , \ |r-q|\leq \gamma\,.\end{eqnarray} 
By (\ref{nesting}) we have
\begin{eqnarray}
\mathcal L_\ve^L (B(\xi_c;r)) \leq \int\psi_{\delta_n/r}\left(\frac{\xi-\xi_c}{r}\right)\mathcal L_\ve^L(d\xi)\,.
\end{eqnarray}
Besides, by (\ref{ballq}) and  by differentiability, there exists $C_n > 0$ such that
\begin{eqnarray}
\psi_{\delta_n/r}\left(\frac{\xi-\xi_c}{r}\right) \leq F^+_{n,k,q}(\xi) + C_n \gamma\,.
\end{eqnarray}
Now, by  (\ref{Nunico}),
\begin{eqnarray}
\lim_\ve \int_{L^2([0,1])} F^+_{n,k,q}(\xi) \mathcal L_\ve^L (d\xi) = \int_{\mathcal C([0,1])}   F^+_{n,k,q}(\xi)\mathbb W(d\xi)
\end{eqnarray}
By (\ref{nesting}) again 
\begin{eqnarray}
\int_{\mathcal C([0,1])}   F^+_{n,k,q}(\xi)\mathbb W(d\xi) \leq \mathbb W (B(g_k, q + \delta_n))\,.
\end{eqnarray}
So far, we have obtained
\begin{eqnarray}
\limsup_\ve \mathcal L_\ve^L (B(\xi_c;r)) \leq \mathbb W (B(g_k, q + \delta_n)) + C_n \gamma\,.
\end{eqnarray}
It remains, in the right hand side, to let $\gamma \to 0$ (hence $g_k \to \xi_c$ and $q\to r$) , and then $n \to \infty$ to get
\begin{eqnarray}
\limsup_\ve \mathcal L_\ve^L (B(\xi_c;r)) \leq \mathbb W (B(\xi_c, r))
\end{eqnarray}
With the same line of reasoning, using the other part of (\ref{nesting}) we can obtain
\begin{eqnarray}
\liminf_\ve \mathcal L_\ve^L (B(\xi_c;r)) \geq \mathbb W (B(\xi_c, r))\,,
\end{eqnarray}
which ends the proof for one ball.

A similar proof can be made for functions approximating intersection of balls as in Theorem 2.2 of \cite{gine1980central} and as a consequence the a.s. weak convergence follows. \qed

Eventually, we have the LDP as in Proposition 3.1. Recall that $(\xi_t^1)$ is the meta-Slepian process defined in (\ref{mS}). We omit the proof since it is the same as in the scalar case.
\begin{prop}
The family $(\mathcal L_\ve)$ satisfies the LDP in $\mathcal M_1(\mathcal C([0,1]))$ equipped with the weak topology,  in the scale $\ve^{-1}$ with good rate function
\begin{eqnarray}
\mathbf{\Lambda}^*(\mathcal L) = \sup_{F \in \mathcal C_b(\mathcal C([0,1]))} \int_{\mathcal C([0,1])} F(\xi)  \mathcal L(d\xi) - \mathbf{\Lambda}(F),
\end{eqnarray}
(the Legendre dual of $\mathbf{\Lambda})$ where for every $F \in \mathcal C_b(\mathcal C([0,1]))$,
\begin{eqnarray}
\mathbf{\Lambda}(F) = \lim_{T\to \infty} T^{-1} \log \mathbb E \int_0^TF(\xi_t^1) dt\,.
\end{eqnarray}
\end{prop}

\section{Discrete versions}
\label{DisV}
For a possible discrete version of Wschebor's theorem and associated LDP, we can consider a continuous process observed at times $(k/n)$ where $k \leq n$ with lag $r$. On this basis, there are two points of view. When 
 $r$ is fixed, there are already results on a.s. convergence of empirical measures of increments of fBm (\cite{BLL-book}) and we explain which LDP holds. 
When $r$ depends on $n$ with $r_n \rightarrow \infty$ and $r_n/n \rightarrow 0$, we are actually changing $t$ in $k/n$   and $\ve$ in $r_n/n$ in the above sections. We state convergence (Prop. \ref{LLNd}) and  LDP (Prop. \ref{LDPd}) 
under specific conditions. 

All the LDPs mentioned take place in $\mathcal M^1(\mathbb R)$ equipped with the weak convergence.

\subsection{Fixed lag}

In \cite{BLL-book}, beyond the Wschbebor's theorem, there are results of a.s. convergence of empirical statistics on the increments of
fBm.  The authors defined p. 39 the second order increments as
\[\Delta_n B_H (i) =  \frac{n^H}{\sigma_{2H}} \left[B_H \left(\frac{i+2}{n}\right) - 2 B_H \left(\frac{i}{n}\right)+ B_H \left(\frac{i}{n}\right)\right].\]
and claimed that as $n\rightarrow \infty$
\begin{eqnarray}
\label{BLL1}
\frac{1}{n-1}\sum_0^{n-2}\delta_{\Delta_n B_H (i)} \Rightarrow \mathcal N \ \ \ (a.s.)\,,
\end{eqnarray}
(Th. 3.1 p.44 in \cite{BLL-book}).
Moreover, in a space-time extension, they proved that 
\begin{eqnarray}
\frac{1}{n-1}\sum_0^{n-2}\delta_{\frac{i}{n}, \Delta_n B_H (i)} \Rightarrow  \lambda\otimes \mathcal N\ \ \ (a.s.)\,,
\end{eqnarray}
(Th. 4.1 in \cite{BLL-EJS}).

Let us restrict for the moment to the case $H=1/2$. The empirical distribution of (\ref{BLL1}) has the same distribution as
\[\frac{1}{n-1}\sum_0^{n-2}\delta_{2^{-1/2}(X_{i+2} - X_{i+1})}\]
where the $X_i$ are independent and $\mathcal N$ distributed.
We can deduce the LDP (in the scale $n$) from the LDP for the 2-empirical measure by contraction. 
If $i$ is the mapping 
\begin{eqnarray}\nonumber\mathbb R^2 &&\rightarrow \mathbb R\\  (x_1,x_2) &&\mapsto (x_2 - x_1)/ \sqrt 2\end{eqnarray}
the rate function is
\begin{eqnarray}
I(\nu) = \inf \{ I_2 (\mu) ;  \mu \circ i^{-1} = \nu \}\,, 
\end{eqnarray}
where $I_2$ is the rate function of the 2-empirical distribution (see \cite{DZ10} Th. 6.5.12).

In the same vein, we could study the LDP for the empirical measure
\[\frac{1}{n-r}\sum_0^{n-r-1} \delta_{\frac{W(k+r)-W(k)}{\sqrt r}}\]
which looks like $\mathcal W_1^r$.
When this lag $r$ is fixed, the scale is $n$ and the rate function is obtained also by contraction ($r=1$ is just Sanov's theorem).

This point of view could be developed also for the fBm using stationarity instead of independence.

\subsection{Unbounded lag}
 Let $(X_i)$ be a sequence of i.i.d. random variables 
 and $(S_i)$ the process of partial sums. Let $(r_n)$ be a sequence of positive integers such that  $\lim_n r_n = \infty$, and assume that
\begin{eqnarray}
\label{C1rn}
\ve_n := \frac{r_n}{n}\searrow 0\,.
\end{eqnarray}
  Set
\begin{eqnarray}
\label{41'}
V_k^n  := \frac{S_{k+r_n} - S_k}{\sqrt{r_n}}\ , \ m_n = \frac{1}{n} \sum_1^n \delta_{V_k^n}\,.
\end{eqnarray}
The next propositions state some extensions of Wschebor's theorem and give the associated LDPs.
The a.s. convergence is obtained only in the Gaussian case under an additional condition.
It seems difficult  to find a general method.
\begin{prop}
\label{LLNd}
\begin{enumerate}
\item If $\mathbb E X_1 = 0, \mathbb E X_1^2 = 1$, then 
\begin{eqnarray}
\label{mnip}
m_n \Rightarrow \mathcal N \ \ \ \hbox{(in probability)}\,.\end{eqnarray}
\item If $X_1 \sim \mathcal N$ and if $(\ve_n)$  
 is 
 such that there exists $\delta \in (0, 1/2)$ and a subsequence $(n_k)$  satisfying
\begin{eqnarray}
\label{assumption}
\sum_k \ve_{n_k} < \infty \ \hbox{and} \ \ve_{n_{k}} = \ve_{n_{k+1}} + o(\ve_{n_{k+1}}^{1+\delta})\,,
\end{eqnarray}
it holds that
\begin{eqnarray}
\label{mnas}
m_n \Rightarrow \mathcal N \ \ \ (a.s.)\,.\end{eqnarray}
\end{enumerate}
\end{prop}
\begin{prop}
\label{LDPd}
\begin{enumerate}
\item Assume that $X_1 \sim \mathcal N$. If  $\lim_n \ve_n n^{1/2} = \infty$, 
then $(m_n)$ satisfies the LDP in the scale 
 $\ve_n^{-1}$ with rate function 
given in (\ref{g15'}-\ref{g15''}) where $\psi = \Psi_1$.
\item Assume that $X_1$ has all its moments finite and satisfies $\mathbb E X_1 = 0$, $\mathbb E X_1^2 = 1$ and that 
 \begin{eqnarray}
\label{liminfsup}
0 < \liminf_n  \ve_n \log n \leq \limsup_n \ve_n \log n < \infty\,.
\end{eqnarray}
Then $(m_n)$ satisfies the LDP in the scale $\ve_n^{-1}$ 
 with rate function
given in (\ref{g15'}-\ref{g15''})  where $\psi = \Psi_1$.
\end{enumerate}
\end{prop}
\begin{rem}
Two  examples of $(r_n)$ satisfying the aassumptions of Prop. \ref{LLNd} 2.  are of interest, particularly in relation to the LDP of Prop. \ref{LDPd}.  The first one is $r_n = \lfloor n^\gamma\rfloor$ with $\gamma \in (0,1)$ (hence $\ve_n \sim n^{\gamma -1}$), for which we can choose $n_k = \lfloor k^{a(1-\gamma)}\rfloor$ with $a > 1$. The second one is $r_n = \lfloor n/\log n\rfloor$ (hence $\ve_n \sim (\log n )^{-1}$), for which we can choose $n_k =\lfloor e^{k^2}\rfloor$. 
\end{rem}

{\bf Proof of Prop. \ref{LLNd}:} We use the method of the above Lemma \ref{az} inspired by \cite{azais1996almost}. For a bounded continuous test function $f$
\[\mathbb E\int f dm_n = \mathbb E f\left(\frac{S_{r_n}}{\sqrt{r_n}}\right) \rightarrow \int fd \mathcal N\]
thanks to the CLT. Moreover
\[\Var\left(\int f dm_n\right)= \frac{1}{n^2} \sum_{|j-k| \leq r_n} \Cov\left(F\left(\frac{S_{j+r_n} - S_j}{\sqrt{r_n}}\right), F\left(\frac{S_{k+r_n} - S_k}{\sqrt{r_n}}\right)\right) \leq \frac{2r_n}{n}||f||_\infty\,.\]
This gives the convergence in probability.

In the Gaussian case, it is possible to repeat the end of the proof of Lemma \ref{az}.
Under our assumption, we see that for any $\beta\in (0, 1/2)$ 
\[\frac{\ve_{n_k}^\beta - \ve_{n_{k+1}}^\beta}{\ve_{n_{k+1}}^{1/2}} =  
 o\left(\ve_{n_{k+1}}^{\delta+ \beta- \frac{1}{2}}\right)\,,\]
which implies that it is enough to choose $\beta \in \left(\frac{1}{2} - \delta, \frac{1}{2}\right)$ . \qed
\medskip

{\bf Proof of Prop. \ref{LDPd}:} 
1) If $X_1 \sim \mathcal N$, then 
\[(V_k^n , k=1, \dots, n) \el \left((\ve_n)^{-1/2}
\left(W\left(\frac{k}{n} + \ve_n 
\right) -W\left(\frac{k}{n}\right)\right)  , k=1, \dots, n\right)\]
and then it is natural to 
 consider  $m_n$ as a Riemann sum. We have now to compare $m_n$ with
\[\mu_{\mathcal W_1^{\ve_n}} = \int_0^1 \delta_{\ve_n^{-1/2}(W(t+ \ve_n) -W(t))} dt\,.\]
It is known that $d_{BL}(\mu, \nu)$ given by (\ref{dBL}) is a convex function of $(\mu, \nu)$ so that : 
\begin{eqnarray*}
&d_{BL} (m_n , \mu_{\mathcal W_1^{\ve_n}}) \leq \displaystyle\int_0^1 d_{BL}(\delta_{\ve_n^{-1/2}(W(t+ \ve_n) -W(t))} , \delta_{V_{\lfloor nt\rfloor}^n}) dt\\
&\leq \ve_n^{-1/2} \displaystyle\int_0^1 \left|W(t+ \ve_n) -W(t) - W\left(\frac{\lfloor nt\rfloor}{n} + \ve_n\right) + W\left(\frac{\lfloor nt\rfloor}{n}\right)\right| dt\\
&\leq 2 (\ve_n)^{-1/2} \sup_{|t-s| \leq 1/n} |W(t) - W(s)|
\end{eqnarray*}
hence
\begin{eqnarray*}
\mathbb P(d_{BL} (m_n , \mu_{\mathcal W_1^{\ve_n}})  > \delta) 
\leq \mathbb P\left(\sup_{|t-s| \leq 1/n} |W(t) - W(s)| > \frac{\delta(\ve_n)^{1/2}}{2}\right) \leq 2 \exp - n \frac{\delta^2\ve_n}{4}\,.
\end{eqnarray*}
If  $\lim_n \ve_n n^{1/2} = \infty$ 
 we conclude that
\[\lim_{n \rightarrow \infty} \ve_n \log \mathbb P (d_{BL}(m_n ,  \mu_{\mathcal W_1^{\ve_n}})> \delta) = -\infty\,,\]
which means that  $(m_n)$ and $(\mu_{\mathcal W_1^{\ve_n}})$ are exponentially equivalent  in the scale $\ve_n^{-1}$ (Def. 4.2.10 in \cite{DZ10}).

Now, 
 from our Prop. 3.1 or 4.1, $(\mu_{\mathcal W_1^{\ve_n}})$ satisfies the LDP in the scale $\ve_n^{-1}$. Consequently, from Th. 4.2.13 of \cite{DZ10}, the family $(m_n)$ satisfies the LDP at the same scale with the same rate function.
\medskip

2) Let us go to the case when $X_1$ is not normal. We use the Skorokhod representation, as in \cite{Heck} or in \cite{March} (see also \cite{CR} Th. 2.1.1 p.88). 

When $(X_i)$ is a sequence of independent (real) random variables such that $\mathbb E X_1=0$ and $\mathbb E X_1^2 = 1$, there exists a probability space supporting a Brownian motion $(B(t); 0 \leq t < \infty)$ and an  increasing sequence $(\tau_i)$ of stopping times
  such that
\begin{itemize}
\item $(\tau_{i+1}- \tau_i)$ are i.i.d., with $\mathbb E \tau_1 = 1$  
\item $(B(\tau_{i+1}) - B(\tau_i))$ are independent and distibuted as $X_1$\,,
\end{itemize}
Moreover, if  $\mathbb E X_1^{2q} < \infty$, then $\mathbb E \tau_1^q < \infty$.

We have \[S_{j+r} - S_j \el B(\tau_{j+r}) - B(\tau_j)\] 
so that 
\begin{eqnarray}m_n \el \tilde m_n :=  \frac{1}{n} \sum_1^n \delta_{\tilde V_k^n} \ \hbox{with} \  \ \tilde V_k^n = \frac{B(\tau_{k+r_n}) - B(\tau_k)}{\sqrt r_n}\,.\end{eqnarray}
We will compare these quantities with
\begin{eqnarray}
\pi_n =  \frac{1}{n} \sum_1^n \delta_{U_k^n} \ \hbox{with} \  \  U_k^n  := \frac{B(k+r_n) - B(k)}{\sqrt{r_n}}\,,\end{eqnarray}
which fall into the regime of the above part of the proof.
We will prove that the sequences $(\tilde m_n)$ and $(\pi_n)$ are exponentially equivalent.

Again by convexity of $d_{BL}$, we have
\begin{eqnarray}
\nonumber
d_{BL}(\tilde m_n, \pi_n) &&\leq \sum_1^n \frac{1}{n} d_{BL}\left(\delta_{\tilde V_k^n}, \delta_{U_k^n}\right)\\
&&\leq \frac{1}{\sqrt{r_n}}\left(\sup_{k\leq n} |B(\tau_{k+r_n}) - B(k+r_n)| + \sup_{k\leq n} |B(\tau_k) - B(k)|\right)
\end{eqnarray}
Our proof will be complete if we show that for all $\delta > 0$
\begin{eqnarray}
\label{CN}
\lim_n 
\frac{r_n}{n} \log \mathbb P\left(\max_{k\leq n + r_n} |B(\tau_k) - B(k)| > \delta \sqrt{r_n}\right) = -\infty\,.
\end{eqnarray}
We will apply three times the following known result.

If $(\xi_i)$ are i.i.d. centered with $\mathbb E(\xi_1)^{2p}) < \infty$ for some $p \geq 1$, then there exists a universal constant $C > 0$ such that for all integers $n \geq 1$
\begin{eqnarray}
\label{ubound}\mathbb E (\xi_1 + \cdots + \xi_n)^{2p} \leq C(2p)! \mathbb E(\xi_1^{2p}) n^p\,,
\end{eqnarray}
(cf.  \cite{March} Lemma 8 or \cite{Heck} Lemma 2.9).

Actually, for $\alpha \in (0,1)$ and $k \leq r_n^\alpha$, with Markov inequality and (\ref{ubound})
\begin{eqnarray}
\label{55}
\mathbb P (|B(\tau_k)| > \delta \sqrt{r_n}) \leq C (2p)! \delta^{-2p} r_n^{-p} 
 \mathbb E((X'_1)^{2p})  k^p
\leq C (2p)! \delta^{-2p} \mathbb E((X'_1)^{2p}) r_n^{(\alpha -1)p}\,,
\end{eqnarray}
and for the same reasons
\begin{eqnarray}
\label{57}
\mathbb P (B(k)| > \delta \sqrt{r_n}) \leq C (2p)! \mathbb E(N^{2p})\delta^{2p} r_n^{(\alpha -1)p}\,.
\end{eqnarray} 
Now, for $k\geq r_n^\alpha$, and $\beta > 1/2$
\begin{eqnarray*}
\mathbb P(|\tau_k - k| \geq k^\beta) \leq  C (2p)!\mathbb E((\tau_1-1)^{2p}) k^{p(1-2\beta)}
\leq C (2p)!\mathbb E((\tau_1-1)^{2p})r_n^{\alpha p(1-2\beta)}\,.
\end{eqnarray*}
Besides,
\begin{eqnarray*}
\nonumber
\mathbb P\left(|B(\tau_k) - B(k)| \geq 2\delta\sqrt{r_n} \ ,\ |\tau_k -k|\leq k^{\beta}\right)  \leq \mathbb P\left(\sup_{|t-k|\leq k^\beta}|B(t) - B(k)| > 2 \delta\sqrt{r_n}\right)\\
\leq 2 \mathbb P\left(\sup_{t \in (0, k^\beta)}|B(t) - B(k)| > 2 \delta\sqrt{r_n}\right)
\leq 4 
e^{ - 2 \delta^2 r_n k^{-\beta}}\,,
\end{eqnarray*}
which, for $k \leq n+ r_n< 2n$, yields
\begin{eqnarray}
\label{61}
\mathbb P\left(|B_{\tau_k} - B_k| \geq 2\delta\sqrt{r_n}, |\tau_k -k|\leq k^{\beta}\right) \leq  4 e^{- 2^{1-\beta} \delta^2 r_n  n^{-\beta}}\,.
\end{eqnarray}
Gathering (\ref{55}-\ref{57}-\ref{61}), we obtain, by the union bound,
\begin{eqnarray}
\nonumber
\mathbb P\left(\max_{k\leq n + r_n} |B(\tau_k) - B(k)| > 2\delta \sqrt{r_n}\right)&\leq 
C_p \left(\delta^{2p} r_n^{1+(\alpha -1)p} + n r_n^{\alpha (1-2\beta)p}\right) \\&+ 8 n e^{- 2^{1-\beta} \delta^2  r_n n^{-\beta}}\,,
\end{eqnarray}
where the constant $C_p > 0$ depends on $p$ and on the distribution of $X'_1$.

Choosing $\beta > 1/2$ and $r_n$ 
 such that
\begin{eqnarray}
\label{C1}
\liminf_n \frac{r_n }{n}\log r_n > 0 \ , \ 
\limsup_n \frac{r_n}{n} \log n < \infty \ , \ 
\liminf_n \frac{r_n^2}{n^{1+\beta}} > 0\,,
\end{eqnarray}
we will ensure that for every $p > 0$ 
\begin{eqnarray} \
\lim \frac{r_n}{n} \log \mathbb P\left(\max_{k\leq n + r_n} |B(\tau_k) - B(k)| > \delta \sqrt{r_n}\right)  \leq - C p
\end{eqnarray}
where $C$ is a constant independent of $p$, which will prove (\ref{CN}).

Now, the set of sufficient conditions  (\ref{C1}) is equivalent to the condition:
\[0 < \liminf_n \frac{r_n}{n} \log n \leq \limsup_n \frac{r_n}{n} \log n < \infty\,,\]
which is exactly (\ref{liminfsup}). 
\qed
\section{Appendix}
The extension of L\'evy's continuity theorem, already invoked in \cite{azais1996almost} is the following. It is probably well known, but since we do not know any reference, we give its proof for the convenience of the reader.
\begin{lem}
Let $\nu_n, \nu$ be probability measures on $\Bbb R^d$ with characateristic functions$\varphi_n, \varphi$. A sufficient condition for $\nu_n \Rightarrow \nu$ is that $\varphi_n({\bf a}) \to \varphi({\bf a})$ for almost every ${\bf a}\in \Bbb R^d$.
\end{lem}

We follow the classical proof as given for instance in Billingsley \cite{billingsleyprobability} Theorem 26.3. Since the  limiting point is determined, we have just to prove the tightness. 
Considering the  compact set $(-M, M]^d \subset\Bbb R^d$, we have by the union bound
\[\nu_n(K^c) \leq \nu_n^k (|x| > M)\]
where $\nu_n^k, k=1, \dots, d$ are the marginals of $\nu_n$. So the problem can be reduced to $d=1$. The basic inequality
\[\nu^k_n (|x] > M) \leq \frac{M}{2}\int_{-2/M}^{2/M} (1 - \nu_n^k(a)) da\,.\]
Since $\varphi$ is continuous at $0$ with $\varphi(0) = 1$, there is for  positive $\eta$ some $M$ such that
\[ \frac{M}{2}\int_{-2/M}^{2/M} (1 - \nu^k(a)) da < \eta\]
Since $\varphi_n$ converges to $\varphi$ almost everywhere, ans since the integrand is bounded by $2$, the dominated convergence theorem gives
\[\lim_n \int_{-2/M}^{2/M} (1 - \nu_n^k(a)) da = \int_{-2/M}^{2/M} (1 - \nu^k(a)) da\,,\]
so that there exists $n_0$ such that 
\[ \frac{M}{2}\int_{-2/M}^{2/M} (1 - \nu_n^k(a)) da < \frac{\eta}{2}\]
for $n\geq n_0$. Ending is routine. \qed

\end{document}